\documentclass{amsart}

\usepackage{graphicx}

\newtheorem{theorem}{Theorem}[section]
\newtheorem{lemma}[theorem]{Lemma}
\newtheorem{corollary}[theorem]{Corollary} 
\newtheorem{proposition}[theorem]{Proposition} 

\theoremstyle{definition}
\newtheorem{definition}[theorem]{Definition}

\newtheorem{fact}[theorem]{Fact} 

\theoremstyle{remark}
\newtheorem{remark}[theorem]{Remark}

\numberwithin{equation}{section}

\newcommand{\R}{\boldsymbol{R}}
\newcommand{\C}{\boldsymbol{C}}
\newcommand{\Z}{\boldsymbol{Z}}
\newcommand\G{{\mathcal G}}
\renewcommand{\phi}{\varphi}
\renewcommand{\epsilon}{\varepsilon}
\newcommand{\op}[1]{{\operatorname{ #1}}}

\begin{document}

\title[Orientability of
linear Weingarten surfaces]{
Orientability
of linear Weingarten surfaces, \\
spacelike CMC-1 surfaces and maximal surfaces
}


\author{Masatoshi Kokubu}
\address{Department of Mathematics,
   School of Engineering,
   Tokyo Denki University\\ 
   Chi\-yo\-da-ku, Tokyo 101-8457,
   Japan}
\curraddr{}
\email{kokubu@cck.dendai.ac.jp}
\thanks{The first author was supported by Grant-in-Aid for 
Scientific Research (C) No.~22540100 
from the Japan Society for the Promotion of Science.
}

\author{Masaaki Umehara}
\address{Department of Mathematical and Computing Sciences,
   Tokyo Institute of Technology
   \\ 
   O-okayama Megro-ku, Tokyo 152-8552,
   Japan}
\curraddr{}
\email{umehara@is.titech.ac.jp}
\thanks{The second author 
was supported by Grant-in-Aid for 
Scientific Research (A) No.~19204005
from the Japan Society for the Promotion of Science.}


    \dedicatory{
Dedicated to Professor Koichi Ogiue on the occasion 
of his seventieth birthday
}

\begin{abstract}
We prove several topological properties of linear 
Weingarten surfaces of Bryant type, 
as wave fronts in hyperbolic $3$-space.
For example, we show the orientability  of
such surfaces, and also co-orientability when 
they are not flat.
Moreover, we show an explicit formula of 
the non-holomorphic hyperbolic Gauss map
via another hyperbolic Gauss map which is holomorphic. 
Using this, we show the orientability and co-orientability
of CMC-1 faces (i.e., constant mean curvature
one surfaces with admissible singular points) 
in de Sitter 3-space.
(CMC-1 faces might not be wave fronts in general, 
but belong to a 
class of linear Weingarten surfaces with singular points.)
Since both linear 
Weingarten fronts and CMC-1 faces may have singular points, 
orientability and co-orientability are both
nontrivial properties. 
Furthermore, we show that the zig-zag representation of
the fundamental group of a linear Weingarten surface 
of Bryant type is trivial.
We also remark on some properties of non-orientable 
maximal surfaces in Lorentz-Minkowski 3-space,  
comparing the corresponding properties of CMC-1 faces
in de Sitter 3-space.
\end{abstract}

\maketitle                   

\section{Introduction}
Let $M^2$ and $N^3$ be 
$C^{\infty}$ manifolds of dimension 2 and of dimension 3, 
respectively.
The projectified cotangent bundle
$P(T^*N^3)$ has a canonical contact structure.
A $C^\infty$ map $f : M^2\to N^3$ is called a {\it frontal} 
if $f$ lifts to 
a Legendrian map $L_f$, i.e.,
a $C^\infty$ map $L_f : M^2 \to P(T^*N^3)$
such that the image $dL_f(TM^2)$ of the tangent bundle $TM^2$ 
lies in the contact hyperplane 
field on $P(T^*N^3)$.
Moreover, $f$ is called a {\it wave front} or a {\it front} if 
it lifts to a Legendrian immersion $L_f$.
Frontals (and therefore fronts) generalize 
immersions, as they allow for
singular points. 
A frontal $f$ is said to be {\it co-orientable}
if its Legendrian lift $L_f$ can lift up to a
$C^{\infty}$ map into the cotangent bundle $T^*N^3$, 
otherwise it is said to be {\it non-co-orientable}.
It should be remarked that, when $N^3$ is Riemannian 
and orientable, 
a front $f$ is co-orientable if and only if 
there is a globally defined unit normal vector field 
$\nu$ along $f$. 
(In \cite{KRUY}, `fronts' were implicitly 
assumed to be co-orientable by definition, and fronts 
which are not necessarily co-orientable 
were distinguished as different, by calling them `p-fronts'. 
However, in this paper, 
a front may be either co-orientable or non-co-orientable and
the term `p-front' is not used.)

In \cite{GMM2}, G\'alvez, Mart\'\i{}nez and Mil\'an
gave a fundamental framework for the theory of
linear Weingarten surfaces of Bryant type
in hyperbolic 3-space $H^3$. 
We shall investigate such surfaces in $H^3$ 
and in de Sitter 3-space $S^3_1$
in the category of (wave) fronts
which admit certain kinds of singular points.
We remark that orientability and co-orientability are not 
equivalent for fronts, though they are equivalent   
for immersed surfaces.
We prove that non-flat linear Weingarten fronts  
of Bryant type in $H^3$ and in  $S^3_1$ are 
co-orientable and orientable. 
Flat fronts, in contrast, are orientable but are not necessarily
co-orientable
(see \cite{KRUY}).
For a co-orientable front $f:M^2\to H^3$,
a representation of the fundamental group
\[
\sigma_f:\pi_1(M^2)\to \Z
\]
called the \lq zig-zag representation\rq\ is induced
(see Section \ref{subsec:zigzag}), which is
invariant under continuous deformations of fronts. 
We shall show  
that
the zig-zag representation of any linear Weingarten front of
Bryant type is trivial (see Theorem \ref{thm:zigzag}).
Moreover, we give an explicit formula for 
the hyperbolic Gauss maps of these surfaces.

On the other hand, we shall also study 
 CMC-1 faces (defined by Fujimori \cite{F}) in $S^3_1$. 
 CMC-1 faces are constant mean curvature one
surfaces on their regular sets, and are
frontals (see Corollary \ref{cor:smooth}), but not fronts in general.
In fact, each limiting tangent plane
at a singular point contains a lightlike vector.
They also belong to a class of linear Weingarten surfaces of
Bryant type and satisfy an Osserman-type
inequality (cf.~\cite{F} and \cite{FRUYY})
under a certain global assumption. 
As an application of an explicit formula for hyperbolic Gauss maps, 
we prove the orientability and the co-orientability
of CMC-1 faces in $S^3_1$,
which is, in fact, deeper than the above corresponding assertion 
for linear Weingarten fronts, because 
CMC-1 fronts in $S^3_1$ admit only isolated singular points 
(see Corollary \ref{cor:s}). 

It should be remarked that
CMC-1 faces in $S^3_1$ have quite similar properties
to maxfaces (i.e. maximal surfaces in Lorentz-Minkowski 3-space
$\R^3_1$
with admissible singular points defined in \cite{UY2}).
In Section 5 (this
section is a joint work with Shoichi Fujimori and Kotaro Yamada),
using the same method as in Section 4, we show
the existence of globally defined real analytic
normal vector field on a given maxface in $\R^3_1$.

\medskip

Throughout this paper, $M^2$ denotes a smooth and
connected  $2$-dimensional manifold or a Riemann surface;  
$\R^4_1$ denotes the Lorentz-Minkowski
space of dimension $4$, with the Lorentz metric
\[
  \bigl\langle (x_0,x_1,x_2,x_3),(y_0,y_1,y_2,y_3)\bigr\rangle =
  -x_0y_0+x_1y_1+x_2y_2+x_3y_3, 
\]
and the hyperbolic $3$-space $H^3$ and the
de Sitter $3$-space $S^3_1$ are defined by 
\begin{align*}
  H^3 &(= H^3_+) = 
  \Bigl\{(x_0,x_1,x_2,x_3)\in\R^4_1 \, ; \, -x_0^2+x_1^2+x_2^2
+x_3^2=-1,
 x_0>0 \Bigr\}, \\
  S^3_1 &=
  \Bigl\{
(x_0,x_1,x_2,x_3)\in\R^4_1 \, ; \, -x_0^2+x_1^2+x_2^2+x_3^2=1
\Bigr\},  
\end{align*}
respectively. 
$H^3$ is a simply-connected Riemannian $3$-manifold with constant
sectional curvature $-1$,  and 
$S^3_1$ is
a simply-connected Lorentzian $3$-manifold with constant
sectional curvature $1$. 
We occasionally consider 
\[
H^3_- := 
  \Bigl\{(x_0,x_1,x_2,x_3)\in\R^4_1 \, ; \, -x_0^2+x_1^2+x_2^2+x_3^2=-1,
 x_0<0 \Bigr\}
\] isometric to $H^3_+$.
\section{Linear Weingarten surfaces of
Bryant type as wave fronts}
Let $f: M^2 \to H^3(\subset \R^4_1)$ be
a front in hyperbolic $3$-space $H^3$.  

If $f$ is co-orientable, then there exists a $C^{\infty}$ map  
$\nu : M^2 \to S^3_1$ satisfying 
$\langle f, \nu \rangle=0$ and
$\langle df, \nu \rangle=0$, 
namely, 
\[
L_f =(f,\nu):M^3 \to T_1H^3 \subset H^3 \times S^3_1
\] 
is the Legendrian lift of $f$. (Here, we identify the unit 
tangent bundle $T_1H^3$ with the unit cotangent bundle $T_1^*H^3$.) 
We call the $C^{\infty}$ map $\nu$ a unit normal field of $f$. 
Interchanging their roles, $\nu$ can be considered 
as a spacelike front in $S^3_1$ 
with $f$ as a unit normal field.
Here, a front in $S^3_1$ is said to be {\it spacelike} if
the plane in $T_{\nu(p)}S^3_1$ orthogonal to the normal
vector $f(p) \in H^3$ is spacelike at each point $p$. 
(Conversely, for a given spacelike co-orientable front in $S^3_1$,
its  unit normal field
is a front which lies in the upper half or the lower half 
of the hyperboloid $H^3_+\cup H^3_-$ of two sheets in 
$\R^4_1$.)  
For each real number $\delta$, the pair 
\begin{eqnarray}\nonumber
f_\delta = (\cosh \delta)f + (\sinh \delta)\nu, \quad
\nu_\delta = (\cosh \delta)\nu + (\sinh \delta)f
\end{eqnarray}
gives a new front $f_{\delta}$ called a {\it parallel front} of $f$
whose unit normal field is $\nu_\delta$.

From now on, we consider a front $f:M^2\to H^3$ which may not
be co-orientable.
The following fact 
implies that even when the front $f$ has a singular point $p$,
taking a co-orientable neighborhood $U$ of $p$,  
we can consider a parallel front 
\[
f_{U,\delta}:=(\cosh \delta)f + (\sinh \delta)\nu
\qquad (\delta\in \R)
\]
on $U$, where $\nu$ is a unit normal vector field
of $f$ on $U$.
For a 
suitable $\delta \in \R$, we may 
assume that 
$f_{U,\delta}$ is an immersion
at $p$. 
\begin{fact}\label{fact:parallel}
Let $f : M^2\to H^3$ be a front. Then for
each $p\in M^2$, taking a co-orientable neighborhood $U$ of $p$, 
the parallel fronts
$f_{U,\delta}$ are immersions at $p$
except for at most two values $\delta\in \R$. 
\end{fact}
First, 
 we begin with the case $f:M^2 \to H^3$ is an immersion, 
which is more restrictive than the case $f$ is a front. 
\begin{definition} \label{def:1}
An immersion $f : M^2\to H^3$ is said to be
 {\it horospherical linear Weingarten} if
the mean curvature $H$ (with respect to the 
normal field $\nu$) and 
the Gaussian curvature $K$
satisfy the relation $a(H-1) + b K=0$ 
for some real constants $a$ and $b$
such that $a^2+b^2\ne 0$. 
\end{definition}
\begin{remark}\label{rem:hlw-im}
A horospherical linear Weingarten immersion $f : M^2 \to H^3$ 
with $a = 0$ is just a flat immersion, hence 
$M^2$ is orientable (see \cite[Theorem B]{KRUY}). On the other hand, 
a horospherical linear Weingarten surface with $a \ne 0$ is 
also orientable by the following reasons: 
\begin{enumerate}
\item[(1)]
If $f$ is isoparametric, then it is orientable because of classification
 of such surfaces. In particular, minimal linear Weingarten 
surfaces orientable since it is isoparametric (see \cite{chen}).
\item[(2)] 
Suppose $f$ is not isoparametric.
 The opposite unit normal field $-\nu$ does not give 
a horospherical linear Weingarten immersion, that is, 
$(-H)-1$ is never proportional to $K$.   
Thus it cannot happen that $M^2$ is non-orientable. 
\end{enumerate}
\end{remark}

\begin{remark}
\label{rem:2}
The parallel surfaces $f_\delta$ of 
a horospherical linear Weingarten
immersion $f$ are also horospherical linear Weingarten
on the set $R_{f_\delta}$ of regular points of $f_\delta$ 
(see \cite{Kok}). 
In particular, parallel surfaces of a flat surface 
are also flat on $R_{f_\delta}$, which is a well known fact. 
\end{remark}
Fact \ref{fact:parallel} and Remark \ref{rem:2} 
enable us to give the following definition.
\begin{definition} \label{def:2}
A front $f : M^2\to H^3$ is called
a {\it horospherical linear Weingarten
front} if for
each $p\in M^2$, there exist a co-orientable neighborhood $U$ of $p$ 
and $\delta_0\in \R$
such that the parallel front
$f_{U,\delta_0}$  is  
a horospherical linear Weingarten immersion on $U$.
(Indeed, the parallel front $f_{U,\delta}$ is a
horospherical linear Weingarten
immersion at $p$
with at most two exceptional values of $\delta \in \R$.)  
\end{definition}

We can give two exceptional examples of 
horospherical linear Weingarten fronts:
\begin{enumerate}
 \item[(1)] a hyperbolic line as a degeneration of
a parallel surface of a hyperbolic cylinder, 
 \item[(2)] a single point as a degeneration of
a parallel surface of a geodesic sphere. 
\end{enumerate} 
These two examples have no regular points 
in their domains of definition. Conversely, 
any horospherical linear Weingarten fronts 
not having regular points are locally
congruent to
one of these two examples
(see the appendix).
So we exclude them from our study in this paper. 
In other words, we will assume (except in the appendix) that 
{\it a horospherical linear Weingarten front has regular points}.

Note that a horospherical linear Weingarten front of zero 
Gaussian curvature 
is  just a flat p-front defined in \cite{KRUY}. 
(See also \cite{KUY} and \cite{KRSUY}.) 
For horospherical linear Weingarten fronts which are not flat, 
the following assertion holds: 

\begin{theorem}\label{thm:non-flat}
Let $f: M^2\to H^3$ be a non-flat horospherical
linear Weingarten front which is not isoparametric. 
Then $f$ is co-orientable.
Moreover, there exists
a unique ratio  $[a:b]$ of real constants 
$(a^2+b^2\ne 0)$
such that the parallel front $f_\delta$ 
satisfies 
\begin{equation}\label{eq:KH}
a(H_\delta-1) + b_\delta K_\delta=0,
\quad b_\delta:=b e^{2 \delta}+\frac{a(e^{2\delta}-1)}2
\end{equation}
on its regular set $R_{f_\delta}$, where $H_\delta$ 
and $K_\delta$ are
the mean curvature and the Gaussian curvature
of $f_\delta$ on $R_{f_\delta}$, respectively.  
\end{theorem}
\begin{proof}
Fix a regular point $p \in M^2$ arbitrarily. 
Let the equation $a(H-1)+bK=0$ ($a^2+b^2 \ne 0$) hold for $f|_U$ 
with the unit normal field $\nu^{(p)}$ 
on a neighborhood $U$ of $p$. 
Then the equation \eqref{eq:KH} holds on $U \cap R_{f_\delta}$ 
because of \cite[(2.18)]{Kok}.

Take another point $q$ close to $p$ so that 
there exist a neighborhood $U_q$ of $q$ and $\delta_0 \in \R$ 
satisfying $U_q \cap U \ne \emptyset$ and 
$f_{U_q,\delta_0}$ is a  
horospherical linear Weingarten immersion with respect to 
the unit normal field $\nu^{(q)}$. 
Then $f_{U_q \cap U,\delta_0}$ is   
horospherical linear Weingarten with respect both to  
$\nu^{(p)}$ and to $\nu^{(q)}$. 
It follows from Remark \ref{rem:hlw-im} that
$\nu^{(q)}$ coincides 
with $\nu^{(p)}$ on $U_q \cap U$. 
Moreover, since the equation \eqref{eq:KH}
holds on $U_q \cap U$, it must hold on $U_q$.

Continuing this argument, we have the collection 
$\{\nu^{(q)} \, | \, q \in M^2 \}$ which then fixes a global 
unit normal field $\nu$. Hence the front $f$ is 
co-orientable.

Finally, we note that the argument above implies  
the equation \eqref{eq:KH} holds on each $U_q \cap R_{f_{\delta}}$, 
therefore on the regular set $R_{f_{\delta}}$.
\end{proof}

\bigskip

Let $f : M^2\to H^3$ be a horospherical 
linear Weingarten front 
which satisfies the relation $a(H-1) + b K=0$
for some real constants $a$ and $b$  ($a^2+b^2\ne 0$). 
Let $\nu : M^2 \to S^3_1$ denote the unit normal field of $f$ 
(assuming the co-orientability whenever $f$ is a flat front). 
Then 
 $f+\nu$ is a map from $M^2$ to the
lightcone $\Lambda^3$ in $\R^4_1$. 
Here, the lightcone $\Lambda^3$ is, by definition, 
\[
  \Lambda^3 =
  \Bigl\{(x_0,x_1,x_2,x_3)\in\R^4_1 \, 
; \, -x_0^2+x_1^2+x_2^2+x_3^2=0 \Bigr\}.
\]
Except for the case where $a+2b=0$, the map $f+\nu$
induces a pseudometric of 
constant Gaussian curvature 
$\varepsilon=a/(a+2b)$ (cf. \cite{GMM2}, \cite{Kok}). 
The equation $a(H-1) + b K=0$ can be rewritten as 
$2\varepsilon (H-1) + (1-\varepsilon) K=0$.

Let $W_\varepsilon(M^2)$ denote the set of 
horospherical linear Weingarten fronts 
from $M^2$ to $H^3$ satisfying
\[
 \begin{cases}
  2\varepsilon (H-1) + (1-\varepsilon) K=0 
& (\text{if $\varepsilon \in \R$}) \\
2 (H-1)- K=0 & (\text{if $\varepsilon = \infty$}).
 \end{cases}
\]

\begin{itemize}
\item (Flat fronts) 
A front $f\in W_0(M^2)$ has zero Gaussian curvature, that
is, $W_0(M^2)$ is the set of flat fronts 
defined on $M^2$
in $H^3$.
(In this case, $a=0$.)
The parallel fronts of a co-orientable flat front $f$ 
also belong to $W_0(M^2)$.
Fundamental properties of flat fronts
are given in \cite{GMM1}, \cite{KRUY} and \cite{KUY}.
The duality between flat surfaces in $H^3$ and those in $S^3_1$
is pointed out in \cite{IS}.
\item (Hyperbolic type)
A front $f\in W_\varepsilon(M^2)$
is called a linear Weingarten front of
{\it hyperbolic type} if $\varepsilon>0$ (i.e., $a/(a+2b)>0$). 
The parallel fronts  
of $f$ are also in the same class 
$\bigcup_{\varepsilon >0} W_{\varepsilon}(M^2)$. 
(cf. \cite{Br}, \cite{UY1}, etc., for the case $\varepsilon=1$.)
For each $f\in \bigcup_{\varepsilon >0} W_{\varepsilon}(M^2)$, 
there exists a unique $\delta \in \R$ such that $f_\delta$ is a CMC-1 
(constant mean curvature one) front
in $H^3$.
In this case, the unit normal vector field
$\nu_{\delta}:M^2\to S^3_1$ is
HMC-1 (harmonic-mean curvature one), that is,
the mean of reciprocals of principal curvatures equals one. 
 (See Corollary \ref{cor:h} for a related result.)

\item (de Sitter type)
A front $f\in W_\varepsilon(M^2)$
is called a linear Weingarten front of
{\it de Sitter type} if $\varepsilon<0$ (i.e., $a/(a+2b)<0$). 
The parallel fronts  
of $f$ are also in the same class 
$\bigcup_{\varepsilon <0} W_{\varepsilon}(M^2)$. 
For each $f\in \bigcup_{\varepsilon <0} W_{\varepsilon}(M^2)$, 
there exists a unique $\delta \in \R$ such that 
$f_{\delta}$ is
HMC-1 (see \cite{Kok}).
In this case, $\nu_\delta$ gives the CMC-1 front in $S^3_1$.
(See Corollary \ref{cor:s}
for a related result.)

\item (Horo-flat fronts) 
A front $f\in W_\infty(M^2)$ 
is said to be {\it horo-flat} (cf. \cite{IST}).
(In this case, $a+2b=0$.)
If $f\in W_\infty(M^2)$, one of 
the hyperbolic Gauss maps 
(see Section \ref{sec:ori}) 
of $f$ degenerates.
The parallel fronts of $f$ also belong to $W_{\infty}(M^2)$.
Fundamental properties of horo-flat fronts are
given in \cite{IST} and \cite{TT}. 
It should be remarked that there exist non-real analytic
horo-flat surfaces (see \cite[Example 2.1]{AG}).
\end{itemize}
\begin{remark}\label{rmk:Iz}
Our terminology \lq {\it horospherical} linear Weingarten\rq\
comes from the following reasons:
The horospheres in $H^3$ satisfy $K=H-1=0$, and 
 are elements in  $W_1(M^2)$, $ W_0(M^2)$,  
$W_{-1}(M^2)$ and $W_\infty(M^2)$ at the same time, as
the degenerate cases.
Recently, Izumiya et al. \cite{IPRT} proposed the horospherical geometry
in $H^3$, which includes all of the above classes $W_\varepsilon(M^2)$. 
\end{remark}
As pointed out at the beginning of this section, 
if $f:M^2\to H^3$ is a co-orientable front,
there exists a unit normal field $\nu$ defined globally on
$M^2$, which gives also a front $\nu: M^2\to S^3_1$.
For $f \in W_{\varepsilon}(M^2)$, one can verify that 
the mean curvature $\hat H$ and the Gaussian curvature
$\hat K$ of $\nu: M^2\to S^3_1$ satisfy
\[
2\varepsilon(\hat H-1)+(1+\varepsilon)\hat K=0.
\]
This implies that flat fronts, horo-flat fronts, 
linear Weingarten fronts of hyperbolic
type and of de Sitter type in $H^3$ 
correspond to those in $S^3_1$, 
respectively. 
\begin{definition}
 Linear Weingarten fronts of hyperbolic type, of de Sitter type 
and flat fronts are all called 
{\it linear Weingarten fronts of Bryant type} 
(cf. \cite{GMM2}). 
In other words, a linear Weingarten front of Bryant type 
is a non-horo-flat linear Weingarten front.
\end{definition}

It should be remarked that Aledo and Espinar \cite{AE} recently
classified complete linear Weingarten immersions of Bryant type
with non-negative Gaussian curvature  in $S^3_1$.
In this paper, almost all formulas are of 
linear Weingarten surfaces in $H^3$, but one can find
the several corresponding formulas 
of linear Weingarten surfaces in 
$S^3_1$ in \cite{AE} (see also \cite{IS}).

\section{Orientability of linear Weingarten fronts 
of Bryant type and an explicit formula for 
the hyperbolic Gauss map $G_*$}\label{sec:ori}
We now concentrate on linear Weingarten fronts 
of Bryant type, i.e., fronts in $W_{\varepsilon}(M^2)$ 
($\varepsilon \ne \infty$), 
which will be denoted by $f : M^2 \to H^3$ 
throughout this section.

\subsection{A representation formula and singular points}

The first, second and third fundamental forms 
are denoted 
by $I$, $I\!I$ and $I\!I\!I$, respectively. 
(In general, $I$ and $I\!I\!I$ are 
well-defined not only for immersions but also for fronts, and 
$I\!I$ is well-defined for co-orientable fronts.)
The sum $I+I\!I\!I$ of the first and third fundamental
forms of $f$ is a positive definite metric on $M^2$, 
because $f$ is a front. 

On the other hand, one can verify that 
the symmetric $2$-tensor $\varepsilon I + (1-\varepsilon) I\!I$ 
for $f \in W_{\varepsilon}(M^2)$ ($\varepsilon \ne 0, \infty$) 
or for a co-orientable $f \in W_0(M^2)$
is definite on $M^2 \setminus S_f$ and vanishes on $S_f$, 
where $S_f$ denotes the set of singular points of $f$. 
We now suppose that $M^2$ is orientable.
Then there is a unique complex structure
on the regular set $R_f:=M^2\setminus S_f$ so that 
$\varepsilon I + (1-\varepsilon) I\!I$ is Hermitian. 
Since this complex structure is common to the whole parallel 
family of $f$,
and the singular set
changes when taking a parallel surface 
(cf. Fact \ref{fact:parallel}),
it follows that
{\it there is a unique complex structure
on $M^2$ so that 
$\varepsilon I + (1-\varepsilon) I\!I$ is Hermitian.}

 $\Lambda^3_+\ni p\mapsto [p]\in \Lambda^3_+/{\sim}$ denotes 
the canonical projection of the positive lightcone 
$\Lambda^3_+ = \{ (x_0,x_1,x_2,x_3) \in \Lambda^3 \, ; \, x_0 >0 \}$ 
onto the set $\Lambda^3_+/{\sim}$ of positively oriented lightlike lines.  
Note that $\Lambda^3_+/{\sim}$ can be naturally identified  
with the Riemann sphere $S^2 \cong \C\cup\{\infty\}$. 
The {\it hyperbolic Gauss map} $G:=[f+\nu]$ can be defined for 
$f \in W_{\varepsilon}(M^2)$ ($\varepsilon \ne 0$) and for a co-orientable 
flat front $f$.  

\medskip

 From now until the end of this subsection,  
we assume that $M^2$ is orientable 
(though we will later know that this assumption is not necessary).  
And also, {\it throughout this section, all flat fronts 
considered are assumed to be co-orientable.}
We have already shown that non-flat linear Weingarten surfaces of
Bryant type are all co-orientable.
On the other hand, even if $f:M^2\to H^3$ 
is a non-co-orientable flat front,
we can take a double covering $\pi:\hat M^2\to M^2$
such that $f\circ \pi$ is co-orientable.
So this assumption of co-orientability is not so restrictive.

\medskip

The hyperbolic Gauss map $G: M^2 \to 
\C\cup\{\infty\}=S^2$ of $f$ is a meromorphic function 
with respect to the complex structure mentioned above.
On the other hand, the metric 
$d\sigma^2$ induced by $f+\nu : M^2\to \Lambda^3_+$ 
is a Hermitian pseudometric on $M^2$
of constant curvature $\varepsilon$
which may have isolated singular points of
integral order (cf. \cite{GMM2}, \cite{Kok}, 
the isolated singular points of $d\sigma^2$ in $M^2$
correspond to umbilical points of $f$).
Then
there exists a meromorphic function
$h$ on $\tilde M^2$ 
satisfying 
\[
d\sigma^2=\frac{4 |dh|^2}{(1+\varepsilon |h|^2)^2},
\]
where $\tilde M^2$ is the universal covering of $M^2$. 
($h$ is the so-called {\it developing map\/}.) 
Then the fundamental forms of $f$ 
are written as follows (cf. \cite{GMM2}, \cite{Kok}): 
\begin{align}\label{eq:I}
I&=
 \frac{(1-\varepsilon)^2}4 d\sigma^2+\frac{4|Q|^2}{d\sigma^2}
+(1-\varepsilon)(Q+\bar{Q}),  \\
\label{eq:II}
I\!I&=
 \frac{\varepsilon^2-1}4 d\sigma^2+\frac{4|Q|^2}{d\sigma^2}
-\varepsilon(Q+\bar{Q}), \\
I\!I\!I &= 
 \frac{(1+\varepsilon)^2}{4}d\sigma^2 
-(1+\varepsilon) (Q + \bar Q)
+ \frac{4 |Q|^2}{d\sigma^2} , \label{eq:third}
\end{align}
where
\begin{equation}
 \label{eq:Q}
Q:= \frac{1}{2} \left(S(h)-S(G)\right)
= \frac{1}{2} \left( \{h:z\}dz^2-\{G:z\}dz^2 \right). 
\end{equation}
The 2-differential $Q$ 
is called the {\it Hopf differential} of $f$.
(The Schwarzian derivative 
$S(G)$ as a locally
defined meromorphic 2-differential
is given by \cite[(2.18)]{KRUY}.
The difference $S(h)-S(G)$ does not depend on 
the choice of local
complex coordinate $z$, see also \cite[Lemma 3.7]{Kok}.)
We also note that 
\begin{align*}
\varepsilon I + (1-\varepsilon) I\!I &= 
-\frac{(1-\varepsilon)^2}4 d\sigma^2+\frac{4|Q|^2}{d\sigma^2},  \\
I + I\!I\!I &= 
\frac{1+\varepsilon^2}{2} d \sigma^2 + \frac{8 |Q|^2}{d\sigma^2} 
-2 \varepsilon ( Q + \bar Q ).  
\end{align*}

The following assertion can be proved 
easily by a modification of 
\cite[Proposition 3.5]{Kok}, which is
a variant of the holomorphic representation
formula in \cite{GMM2}:

\begin{theorem}\label{thm:kok-rep}
Let $G$ be a meromorphic function on a Riemann surface $M^2$,
and $h$ a meromorphic function on the universal covering $\tilde M^2$
such that the pseudometric $d\sigma^2:=4|dh|^2/(1+\varepsilon |h|^2)^2$
gives a single-valued symmetric tensor on $M^2$.
Suppose that 
\begin{equation}\label{eq:front-condition}
\frac{1+\varepsilon^2}{2} d \sigma^2 + \frac{8 |Q|^2}{d\sigma^2} 
-2 \varepsilon ( Q + \bar Q ) 
\end{equation}
is positive definite, where 
$Q:= (S(h)-S(G))/2 = (\{h:z\}dz^2-\{G:z\}dz^2)/2$. 
Define $f=\G \mathcal A \mathcal \G^*$ and
$\nu=\G \mathcal B \mathcal \G^*$ by
\[
\G:=i(G_h)^{-{3}/{2}}
\left(
\begin{array}{cc}
-GG_h & (GG_{hh}/2)-(G_h)^2 \\
-G_h & G_{hh}/2
\end{array}
\right), 
\]
where $G_h:=dG/dh,G_{hh}:=dG_h/dh$
and
\begin{eqnarray*}
\mathcal A:=
\left(
\begin{array}{cc}
 \frac{1+\varepsilon^2 |h|^2}{1+\varepsilon |h|^2} & 
-\varepsilon \overline{h} \\
-\varepsilon h & 1+\varepsilon |h|^2
\end{array}
\right), \quad 
\mathcal B:=
\left(
\begin{array}{cc}
\frac{1-\varepsilon^2 |h|^2}{1+\varepsilon |h|^2} & 
\varepsilon \overline{h} \\
\varepsilon h & -1-\varepsilon |h|^2
\end{array}
\right).
\end{eqnarray*}
Then $f : M^2 \to H^3 = SL(2,\C)/SU(2)$ 
is a co-orientable
linear Weingarten front of Bryant type in $W_\varepsilon(M^2)$, 
and $\nu : M^2 \to S^3_1= SL(2,\C)/SU(1,1)$ is
its unit normal. 
Moreover, $G$ and $($\ref{eq:front-condition}$)$ coincide 
with the hyperbolic Gauss map
and the sum $I + I\!I\!I$
of the first and third fundamental forms of $f$,
respectively.
 Conversely, 
any co-orientable linear Weingarten front of Bryant type in 
$W_{\varepsilon}(M^2)$, 
except for the horosphere, is given in this manner.
\end{theorem}
\begin{remark}\label{rmk:unique}
We give some additional explanation for Theorem \ref{thm:kok-rep}:
\begin{enumerate}
\item
We often identify $\R^4_1$ with the set $\mathrm{Herm}(2)$ of 
$2\times 2$ Hermitian matrices
by
\begin{equation}\label{eq:herm-mink}
   X=(x_0,x_1,x_2,x_3)\leftrightarrow
    X=\sum_{k=0}^3 x_k e_k
    =\begin{pmatrix}
      x_0+x_3           & x_1+i x_2 \\
      x_1-i x_2 & x_0-x_3
     \end{pmatrix},
\end{equation}
where $i=\sqrt{-1}$ and 
\[
    e_0=\begin{pmatrix}1&0\\0&1\end{pmatrix},~
    e_1=\begin{pmatrix}0&1\\1&0\end{pmatrix},~
    e_2=\begin{pmatrix}\hphantom{-}0&i\\-i&0\end{pmatrix},~
    e_3=\begin{pmatrix}1&\hphantom{-}0\\0&-1\end{pmatrix}. 
\]
Since the Lorentzian inner product $\langle \ , \ \rangle$ is 
given by 
\[
   \langle X,Y\rangle
    =-\frac{1}{2}\operatorname{trace}\left(Xe_2 {}^t Y e_2\right),
    \qquad
   \langle X,X\rangle=-\det X, 
\]
 $H^3$ and $S^3_1$ are given by 
\begin{align*}
 H^3 &= \{X\,;\,X^*=X\,,\det X=1, \mathrm{trace} X >0 \}, \\
 S^3_1 &= \{X\,;\,X^*=X\,,\det X=-1\},              
\end{align*}
respectively, where $X^*$ denotes the transposed conjugate 
${}^t \! \bar{X}$ 
to the matrix $X$. 
 Hence, they are also represented by  
\begin{align*}
  H^3 &= \{aa^*\,;\,a \in SL(2,\C) \} = SL(2,\C)/SU(2), \\
 S^3_1 &= \{ae_3a^*\,;\,a\in SL(2,\C) \} = SL(2,\C)/SU(1,1).  
\end{align*}
 \item The choice of the developing map $h$
of $d\sigma^2$ is not unique.  
However, the above expression of $f$ does not depend on
the choice of $h$.
We prove this here:

Suppose that we have the following two expressions
\[
d\sigma^2=\frac{4 |dh|^2}{(1+\varepsilon |h|^2)^2}
=\frac{4 |d\tilde h|^2}{(1+\varepsilon |\tilde h|^2)^2}.
\]
By taking a parallel front, we may assume that $f$
is an immersion on an open subset $U$ of $M^2$.
Since $f$ is real analytic, it is sufficient to show that
$\tilde f$ associated to the pair $(G,\tilde h)$
coincides with  $f$.
Since the first and second fundamental forms $I$
and $I\!I$ can be written in terms
of $G,\,\, d\sigma^2$ and $Q$ as in \eqref{eq:I} and \eqref{eq:II},  
$\tilde f$ is congruent to $f$ by the fundamental theorem
of surface theory.
Moreover, we can show that $f=\tilde f$ in this situation.
In fact, $a=(a_{ij})_{i,j=1,2}\in \op{SL}(2,\C)$ isometrically 
acts $f$ as $afa^*$. Then the hyperbolic Gauss map
$G$ changed as $(a_{11}G+a_{12})/(a_{21}G+a_{22})$.
Since 
$\tilde f$ and $f$ have the same hyperbolic Gauss map $G$,
which implies that $f=\tilde f$.
\item
$\mathcal{G}$ satisfies 
\begin{equation}\label{eq:F-1dF}
 \mathcal{G}^{-1} d \mathcal{G} = 
\begin{pmatrix}
 0 & \theta \\ dh & 0
\end{pmatrix},  
\end{equation}
where $\theta = - \dfrac{1}{2} \{ G : h \} dh $. 
Moreover, note that $Q = \theta dh$. (See \cite{GMM2} and \cite{Kok}
for details.)
\end{enumerate}
\end{remark}

\begin{corollary}\label{cor:h}
Let $f\in W_\varepsilon(M^2)$ for $\varepsilon>0$.
Then 
there exists a unique real number
$\delta$ such that $f_\delta : M^2\to H^3$ is a CMC-1
front 
whose singular set $S_{f_\delta}$ 
consists only of isolated points in $M^2$. 
\end{corollary}
\begin{proof}
Without loss of generality, we  may assume that
$f$ is not a horosphere.
There exists a unique real number
$\delta$ such that $f_\delta\in W_{1}(M^2)$. 
It follows from \eqref{eq:I} with $\varepsilon=1$
that the first fundamental form of $f_\delta$ is
$
I_{f_\delta}={4|Q|^2}/{d\sigma^2},
$
which is positive definite except at the 
zeros of the holomorphic 2-differential $Q$.
Since $Q$ is holomorphic, 
 $f_\delta$ is a front with only isolated singular points.
\end{proof}

Similarly, we get the following:

\begin{corollary}\label{cor:s}
Let $f\in W_\varepsilon(M^2)$ for $\varepsilon<0$.
Then there exists a unique real number
$\delta$ such that the unit normal field $\nu_\delta$
of $f_\delta$ gives a CMC-1 front $\nu_\delta : M^2\to S^3_1$
whose singular set $S_{\nu_\delta}$ 
consists only of isolated points in $M^2$. 
\end{corollary}
\begin{proof}
Without loss of generality, we  may assume that
$f$ is not a horosphere.
There exists a unique real number
$\delta$ such that $f_\delta\in W_{-1}(M^2)$, that is, 
$f_{\delta}$ has constant harmonic-mean 
curvature one (cf. \cite{Kok}). At the same time,  
$\nu_\delta : M^2\to S^3_1$ is a CMC-1 front. 
It follows from \eqref{eq:I} and \eqref{eq:third}
with $\varepsilon=-1$ that  
the first fundamental
form of $\nu_\delta$ is
$
I_{\nu_\delta}=I\!I\!I_{f_\delta}=
{4|Q|^2}/{d\sigma^2}. 
$
Since $Q$ is holomorphic, 
$\nu_\delta$ is a front with only isolated singular points.
\end{proof}

\begin{remark}\label{rmk:isolated}
A meromorphic function $G=z+iz^2$ and 
pseudometrics 
\[
d \sigma^2 = \frac{4|dh|^2}{(1+ |h|^2)^2} \quad
(\mbox{resp.} \,\, d \sigma^2 = \frac{4|dh|^2}{(1- |h|^2)^2}),
\] 
where $h=z+z^3$, induce CMC-1 fronts in $H^3$ 
(resp. in $S^3_1$). 
They have an isolated singular point at $z=0$ because 
 $S(h)-S(G)(=2Q)$ vanishes at $z=0$.
\end{remark}

\begin{remark}\label{rmk:isolated2}
Let $f\in W_0(M^2)$ and $\nu$ its unit normal field.
Then one can show that
$f+\nu$ is a flat front in the lightcone
in $\R^4_1$,
whose singular set $S_{f+\nu}$ 
consists only of isolated points. 
\end{remark}

More generally, we can show the following propositions 
concerning singular points of linear Weingarten 
fronts of Bryant type.
\begin{proposition}
  Let $f : M^2 \to H^3$ be a linear 
Weingarten front of Bryant type written in terms of 
$(G, d \sigma^2=4|dh|^2/(1+\varepsilon|h|^2)^2)$ 
as in Theorem 
\ref{thm:kok-rep}.  
Then the set of singular points
equals 
\[
S_f=  \biggl\{ p \in M^2 \, ; \,
\varepsilon I +(1-\varepsilon) I\!I =0
\biggr\}
=
 \biggl\{ p \in M^2 \, ; \, \frac{4 |Q|^2}{d \sigma^2} -  
\frac{(1-\varepsilon)^2}{4} d\sigma^2=0 \biggr\}, 
\]
where 
$Q$ is the Hopf differential \eqref{eq:Q}. 
Moreover, 
a singular point $p \in S_f$ is nondegenerate if and only if
\begin{enumerate}
 \item[(1)]  $f$ is not CMC-1, and 
 \item[(2)] 
$4\varepsilon h_z \bar{h}  
+(1+\varepsilon |h|^2)(\hat \theta_z / \hat \theta  - h_{zz}/h_z) \ne 0$ 
holds at $p$,  
where 
$ \theta = \hat \theta dz$  
 for a local coordinate $z$, and $h_z =dh/dz$, 
$h_{zz} = d^2 h/dz^2$. 
\end{enumerate} 
\end{proposition}

We now set
\[
\varDelta := \operatorname{Im} \left[ 
\dfrac{1}{\sqrt{1 -\varepsilon}}
\left\{ \dfrac{4 \varepsilon h_z \bar h}{1+\varepsilon |h|^2} + 
\dfrac{\hat \theta_z}{\hat \theta} - \dfrac{h_{zz}}{h_z} \right\}
\dfrac{1}{\sqrt{h_z \hat \theta}}
\right],
\]
where $\sqrt{1 -\varepsilon}$ is an imaginary number when 
$1-\varepsilon<0$.

\begin{proposition}
 For a nondegenerate singular point $p$ of a 
linear Weingarten front $f : M^2 \to H^3$ of Bryant type, 
the germ of $f$ at $p$ is locally diffeomorphic to 
\begin{enumerate}
 \item[(1)] a cuspidal edge if and only if $\varDelta \ne 0$ at $p$, or
 \item[(2)] a swallowtail if and only if $\varDelta = 0$ at $p$ and 
$\dfrac{d}{dt} \Bigr|_{t=0} \varDelta \circ \gamma \ne 0$, where
$\gamma$ is a parametrization of $S_f$ with $\gamma(0) = p$. 
\end{enumerate}
\end{proposition}
The proofs of the above propositions for the case of $\varepsilon =0$ 
are given in \cite{KRSUY}. The proofs for the case of $\varepsilon \ne 0$ 
can be followed by a quite similar argument to that in \cite{KRSUY},
 though we need very lengthy calculations. We omit them here.

\subsection{Zig-zag representation and orientability}\label{subsec:zigzag}
For a given co-orientable front $f : M^2 \to H^3$, 
the {\it zig-zag representation}
\[
\sigma_f : \pi_1(M^2)\to \Z
\]
is induced (cf. \cite{LLR}, \cite{SUY1} 
and \cite{SUY2}), which is invariant under 
the deformation of $f$ as a wave front, 
and 
\begin{equation}\label{eq:equal}
\sigma_f=\sigma_\nu
\end{equation} 
holds for 
the unit normal field $\nu : M^2 \to S^3_1$.

\begin{theorem}\label{thm:zigzag}
Let $f : M^2\to H^3$ $($resp. $f : M^2\to S^3_1)$
be a co-orientable
linear Weingarten front of
Bryant type. Then $M^2$ is orientable.
Moreover, the zig-zag representation $\sigma_f$ is trivial.
\end{theorem}

\begin{figure}[hb]
\begin{center}
        \includegraphics[height=1.5cm]{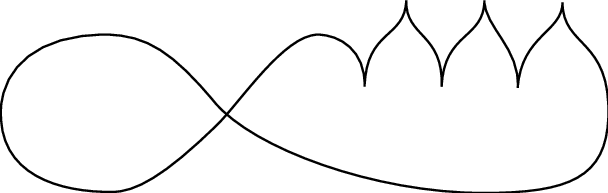}\hspace{1cm}
\caption{A planar curve of zig-zag number 3}\label{zig}
\end{center}
\end{figure}

\begin{proof}
It is sufficient to prove the assertion for
a co-orientable linear Weingarten front
$f : M^2\to H^3$.
(In fact, suppose that the assertion holds for 
linear Weingarten front of
Bryant type in $H^3$.
Let $f : M^2\to S^3_1$ be 
a co-orientable linear Weingarten front.
Then its unit normal vector field 
$\nu:M^2\to H^3$ is 
a co-orientable
linear Weingarten front of
Bryant type.
Thus $M^2$ must be orientable. 
On the other hand, the zig-zag representation is
trivial by \eqref{eq:equal}.)

Suppose that $f:M^2 \to H^3$ is of hyperbolic type (resp.
de Sitter type). 
By Corollary \ref{cor:h} (resp. Corollary \ref{cor:s}),
there exists $\delta\in \R$
such that  $f_\delta$ 
(resp. $\nu_\delta$) is a CMC-1 front with only isolated singular points.
Since $f_\delta$ (resp. $\nu_\delta$)
is an immersion of non-vanishing mean curvature
on $R_{f_\delta}=M^2\setminus S_{f_\delta}$, the open 
submanifold $R_{f_\delta}$ is orientable. 
Since $S_{f_\delta}$ is discrete, we conclude that  
$M^2$ is also orientable. 
On the other hand, we can take a loop $\gamma$ which 
represents a given element of $\pi_1(M^2)$ so that
$\gamma$ does not pass through $S_{f_\delta}$.
Since the zig-zag representation is trivial when the 
corresponding loop does not meet the 
singular set,
$\sigma_{f_\delta}([\gamma])=\sigma_{\nu_\delta}([\gamma])=\sigma_f([\gamma])$
vanishes.

Next, we consider the case 
where $f$ is a co-orientable flat front. 
The orientability of $M^2$ has already been proved in \cite{KRUY}.
We set (see \S3 of \cite{KRSUY})
\[
|\rho_{\delta}|:=
e^{-2\delta}|\theta/dh|=e^{-2\delta}|Q/dh^2|. 
\]
Then the singular set of the parallel
front $f_\delta$ is given by the
$1$-level set $|\rho_{\delta}|=1$.
We fix an element $[\gamma]\in \pi_1(M^2)$,
where $\gamma$ is a loop in $M^2$.
Then we can choose $\gamma$ so that it does not pass through
the zeros of $dh$ nor of $\theta$.
Then there exists a constant $c>0$ such that
$|Q/dh^2|\ge c$ on $\gamma$. 
If we take $\delta$ sufficiently small,
then  we have $e^{-2\delta}|Q/dh^2|>1$ on $\gamma$.
Then $f_\delta$ has no singular point on $\gamma$.
Thus $\sigma_{f_\delta}([\gamma])=\sigma_f([\gamma])$ vanishes.
\end{proof}
\begin{remark}\label{rem:existence}
In contrast to the case of $H^3$, 
there are many flat fronts in $\R^3$
with non-trivial zig-zag representations. 
For example, consider a cylinder over the planar curve 
in Figure \ref{zig}.
\end{remark}
\begin{remark}\label{rem:minimal}\rm
An immersion $f:M^2\to \R^3$ is called a {\it  
linear Weingarten
surface of minimal type}
if the mean curvature $H$ and the Gaussian curvature $K$
satisfy the relation $H=aK$, where $a$ is a real constant.
Such surfaces in $\R^3$ are the corresponding analogue of
linear Weingarten surfaces of Bryant type.
As in the method of Section 1, one can define 
{\it linear Weingarten fronts of minimal type}.
One can easily prove that each
linear Weingarten front $f:M^2\to \R^3$ of minimal type
contains a minimal surface $f_0:M^2\to \R^3$ as a 
wave front in its parallel family. 
Since minimal surfaces have
a well-known Weierstrass representation formula,
the first fundamental form of $f_0$ is given by
$|Q|^2/d\sigma^2$, where $Q$ is the Hopf differential and 
$
d\sigma^2:={4|dg|^2}/(1+|g|^2)^2
$
for the Gauss map $g$ of $f$ (and $f_0$).
Since the singular points of $f_0$ are isolated,
the orientability of $M^2$ is equivalent to the
co-orientability of $f$.
(This is a different phenomenon from the case of
linear Weingarten fronts of Bryant type.
In fact, non-orientable minimal immersions are known.)
By the same argument as in the above proof,
the zig-zag representation of 
a given  linear Weingarten front of minimal type
$f:M^2\to \R^3$ vanishes if $f$ is
co-orientable.
\end{remark}

\subsection{An explicit formula for $G_*$}
Recall that the hyperbolic Gauss map of a 
co-orientable front $f$ is given by $G=[f+\nu]$.
On the other hand, we set $G_*:=[f-\nu]$, which is also
called  the  {\it $($opposite$)$ hyperbolic Gauss map}.
For a co-orientable flat front, it is known that $G_*$ is also
a meromorphic function, 
as well as $G$ (cf. \cite{GMM1}).  We give here
an explicit formula for $G_*$: 

\begin{proposition}\label{prop:Gstar} 
Under the notation of Theorem \ref{thm:kok-rep}, 
the opposite hyperbolic Gauss map $G_*$ is given by 
\begin{equation}\label{eq:G_h}
G_*=G-\frac{(G_h)^2(1+\varepsilon|h|^2)}{\varepsilon \bar h G_h +(G_{hh}/2)
(1+\varepsilon |h|^2)}. 
\end{equation}
\end{proposition}

\begin{proof}
The matrices  $\mathcal A$ and $\mathcal B$ 
in Theorem \ref{thm:kok-rep}
split as 
\[
 \mathcal A = \Phi \Phi^*, \quad 
\mathcal B= \Phi e_3 \Phi^*,
\]
where $\Phi^*:={}^t  \bar{\Phi}$ and
\[
\Phi:= \frac{i}{\sqrt{1 + \varepsilon|h|^2}}
\left(
\begin{array}{cc}
  -1 & - \varepsilon \bar h \\ 0 & 1+ \varepsilon |h|^2
\end{array}
\right), \quad 
e_3:=
\left(
\begin{array}{cc}
    1 & 0 \\ 0 & -1 
\end{array}
\right).
\]
Hence we obtain 
\[
 f - \nu = \G (\mathcal A - \mathcal B) \G^* 
= 2 \, \G \Phi 
\left(
\begin{array}{cc}
    0 & 0 \\ 0 & 1
\end{array}
\right)
(\G \Phi)^*. 
\]
Setting $\G \Phi = 
\left(
\begin{array}{cc}
	 p & q \\ r & s
\end{array}
\right)$, we have  
$f -\nu = 2 
\left(
\begin{array}{c}
 	  q \\ s 
\end{array}
\right)
\left(
\begin{array}{cc}
 	\bar q & \bar s 
\end{array}
\right)$. 
This implies 
$
 G_* = [f -\nu] = q/s \in \C \cup \{ \infty \}.
$ 
Then by a straightforward calculation, we get
the assertion.
\end{proof}
We immediately get the following:
\begin{corollary}\label{cor:G_*}
$G_*$ is meromorphic, as well as $G$, 
if and only if $\varepsilon=0$.
\end{corollary}

Let $D$ be the diagonal set of $S^2\times S^2$.
Then the space $L(H^3)$ 
of oriented geodesics of hyperbolic 3-space $H^3$
can be identified with $S^2\times S^2\setminus D$ 
which is a parahermitian symmetric space 
$SO^0(1,3)/(SO(2)\cdot \R^*)$
 (see Kanai \cite{Kanai}, Kaneyuki \cite{Kane}
and also \cite{CFG}). 
Recently, Georgiou and Guilfoyle 
\cite{GG} proved that 
$L(H^3)$ has a canonical neutral K\"ahler structure
and the pair of hyperbolic Gauss maps 
\[
(G,G_*) : M^2\to S^2\times S^2\setminus D
\]
gives a Lagrangian surface with zero Gaussian curvature
if $f$ is a Weingarten surface.
Since our explicit formula \eqref{eq:G_h} for $G_*$ is written 
in terms of only $G$, $h$ and their derivatives, 
it might be useful for constructing 
such surfaces with positive genus.

\section{Orientability and co-orientability of CMC-1 faces in $S^3_1$}

In this section, we shall give an application of
our explicit formula (\ref{eq:G_h}) for $G_*$.
Let $M^2$ be a Riemann surface.
A holomorphic map $F:M^2\to SL (2,\C)$ is called
{\it null} if the determinant of the derivative
$dF/dz$ with respect to each local complex coordinate $z$
vanishes identically.
It is well known that
CMC-1 surfaces in $H^3$
and in $S^3_1$ are both  projections of null holomorphic curves 
 $F: \tilde{M^2} \to SL (2,\C)$,  
where $\tilde M^2$ is the universal covering of $M^2$ and
$ H^3=SL(2,\C)/SU(2)$ and 
$S^3_1=SL(2,\C)/SU(1,1)$. 
Suppose that a CMC-1 surface in $S^3_1$ is given by 
$f= -\mathcal{GBG}^*$ for $\varepsilon =-1$
as in Theorem \ref{thm:kok-rep}. 
(Here we suppose $f= -\mathcal{GBG}^*$ instead of $\mathcal{GBG}^*$ 
because it makes the argument below more compatible 
to \cite{F} and \cite{FRUYY}. $-f$ and 
$f$ differ merely by congruence in $S^3_1$, that is, 
they are essentially the same.)
Then the null holomorphic lift $F$ 
of $f$ is given by 
\begin{equation} \label{eq:null-lift}
 F = \mathcal{G} \begin{pmatrix}
		  0 & -i \\ -i & ih
		 \end{pmatrix}
\end{equation}
Conversely, the projection 
$Fe_3 F^* : M^2\to S^3_1$ may admit singular points. 
Such a surface is called a CMC-1 face. More precisely, 
a $C^{\infty}$ map $f : M^2\to S^3_1$ is called a {\it CMC-1 face}
if the regular set is open dense in $M^2$, and
for each $p\in M^2$, there exist a neighborhood $U$
of $p$ and a null holomorphic immersion $F : U\to SL(2,\C)$
such that it has an expression $f=Fe_3 F^*$ on $U$ 
(see \cite{F} and \cite{FRUYY}). 
Note that we usually take a null holomorphic lift $F$ 
defined on 
the universal cover $\tilde{M}^2$ for a CMC-1 face. 

We shall show in this section the orientability
and the co-orientability of CMC-1 faces in this setting.
To show the orientability, we will need the explicit formula
for $G_*$, which was a formula missing  in
\cite{FRUYY}.
It should also be remarked that CMC-1 faces 
are not fronts in general but are all frontals,  
as seen below (see \eqref{eq:fuji1}).  
On the other hand, a CMC-1 front in $S^3_1$ may admit
isolated singular points in $M^2$, as seen in Corollary \ref{cor:s},
although CMC-1 faces do not.
This implies that the set of CMC-1 fronts is
not included in the set of CMC-1 faces and vice versa.

For a CMC-1 face $f=Fe_3F^*(=-\mathcal{GBG^*})$, 
we have the following:
\begin{eqnarray}\label{eq:F1}
F^{-1}dF=
\left(
\begin{array}{cc}
h & -h^2 \\ 1 & -h
\end{array}
\right)
\dfrac{Q}{dh}, \quad 
dFF^{-1}=
\left(
\begin{array}{cc}
G & -G^2 \\ 1 & -G
\end{array}
\right)
\dfrac{Q}{dG},
\end{eqnarray}
where $G$, $h$ and $Q$ 
are corresponding to the data in Theorem \ref{thm:kok-rep}. 
The formulas \eqref{eq:F1} follow from 
\eqref{eq:F-1dF} and \eqref{eq:null-lift}. 
Note that the Hopf differential $Q$ and the secondary Gauss map $g$   
have already been introduced in \cite{F}. 
In our notation, the Hopf differential is denoted by the same $Q$, 
the secondary Gauss map $g$ coincides with $h$,
and $\omega:=Q/dg$ in \cite{F} is equal to $\theta =Q/dh$. 

The singular set $S_f$ of a CMC-1 face 
$f=Fe_3F=-\mathcal{GBG^*}$ 
is given by 
\[
 S_f = \left\{ p \in M^2 \, ; \, |h(p)|=1 \right\}.
\]
(See \cite[Theorem 1.9]{F}.) Away from the singular set $S_f$, 
we can define a unit normal field 
$\nu:{M^2\setminus S_f}\to H^3_+\cup H^3_-$ by
$\nu=\mathcal{G}\mathcal A \mathcal{G}^*$. 
It follows from \eqref{eq:null-lift} that 
\begin{equation}\label{eq:fuji}
\nu=
F 
\left(
\begin{array}{cc}
 i h & i \\
 i & 0
\end{array}
\right)
\mathcal A 
{\left(
\begin{array}{cc}
 i h & i \\
 i & 0
\end{array}
\right)}^{\! \! *} F^*
=\frac{\tilde \nu}{1-|h|^2},
\end{equation}
where
\begin{equation}\label{eq:fuji1}
\tilde \nu:=F
\left(
\begin{array}{cc}
 1+|h|^2 & 2h \\ 2 \bar h & 1+|h|^2
\end{array}
\right)
F^*.
\end{equation}
(See also \cite[Remark 1.2]{F}.)
Since the equivalence class $[\tilde \nu]$ 
gives a section of $M^2$ 
into the projective tangent 
bundle $P(TS^3_1)$,   
an arbitrary CMC-1 face $f$ is a frontal.
(The definition of a frontal is given in the introduction.)
Generic singular points on CMC-1 faces are cuspidal edges,
swallowtails and cuspidal cross caps (see \cite{FSUY}). 
 However, cuspidal cross caps
never appear on fronts. Hence, CMC-1 faces are not fronts in general.

\begin{theorem}\label{thm:cmc-face}
Let $f : M^2\to S^3_1$ be a CMC-1 face.
Then  $M^2$ is orientable.
\end{theorem}

\begin{proof}
Suppose, by way of contradiction, 
that $M^2$ is not orientable. 
There exists a double covering 
$\hat\pi : \hat M^2 \to M^2$ such that
$\hat M^2$ is orientable.
Set $\hat f:=f\circ \hat \pi$. 
Let $T : \hat M^2\to \hat M^2$ denote 
the covering involution which is orientation-reversing.
We fix a regular point $p\in \hat M^2$. 
Then there exists a simply connected
neighborhood $U(\subset \hat M^2)$ of $p$ such that
the restriction $\hat f|_U$ is an immersion.
Let $\hat \nu$ be the unit normal 
field on $U$, 
and $\hat G$, $\hat G_*$ denote the hyperbolic Gauss maps 
for $\hat f$ on $U$. 
Then $\hat f\circ T=\hat f$, and either
$\hat \nu \circ T(q)=\hat \nu(q)$ or 
$\hat \nu \circ T(q)=-\hat \nu(q)$
holds for each $q\in U$ as vectors in $\R^4_1$.
First, we suppose 
$\hat \nu\circ T=-\hat \nu$.
Then
\[
\hat G\circ T=[\hat f\circ T+\hat \nu\circ T]
=[\hat f-\hat \nu]=\hat G_*
\]
holds on $U$.
This implies that $\hat G_*$ is 
anti-holomorphic, which contradicts
Corollary \ref{cor:G_*}.
Next, we suppose $\hat \nu\circ T=\hat \nu$.
Then, we have
\[
\hat G\circ T= [ \hat f\circ T+\hat \nu\circ T] =
[ \hat f+\hat \nu ] =\hat G,
\]
which contradicts that $T$ is orientation-reversing. 
\end{proof}

Since $G$ is a holomorphic map,
when the unit normal field
$\nu$ of a CMC-1 face $f$ 
crosses a singular curve $S_f = \{ |h|=1 \}$, 
the image of $\nu$ moves 
into another sheet of the hyperboloid $H^3_+\cup H^3_-$.
Thus, it is natural to expect $\nu$
to be smooth at the singular set under a certain
compactification of 
$H^3_+\cup H^3_-$.   
The {\it hyperbolic $3$-sphere} $\bar {\mathcal H}^3$ is 
a $3$-dimensional manifold 
diffeomorphic to the 3-sphere
\[
 \bar {\mathcal H}^3:=\R^3\cup\{\infty\} \cong S^3
\]
endowed with the metric
$4dx\cdot dx /(1-|x|^2)^2$ 
on $S^3 \setminus \{ {\rm equator} \}$, 
where $x:=(x_1,x_2,x_3) \in \R^3\cup\{\infty\}$.
We consider the stereographic projection
\begin{equation}\label{eq:phi}
 \varphi:H^3_+\cup H^3_- \ni (x_0,x_1,x_2,x_3)\mapsto 
\frac{(x_1,x_2,x_3)}{1-x_0}\in 
\bar {\mathcal H}^3,
\end{equation}
which is an isometric embedding, 
and $\bar {\mathcal H}^3$ can be considered
as a compactification of $H^3_+\cup H^3_-$.

\begin{theorem}\label{thm:smooth}
The unit normal field 
of a CMC-1 face can be considered as a real analytic map
into the hyperbolic $3$-sphere $\bar {\mathcal H}^3$.
In other words, the map $\varphi \circ \nu$, 
which is defined on $M^2 \setminus S_f$,
can be real analytically extended across the singular set $S_f$. 
\end{theorem}

\begin{proof}
We set $F=
\left(
\begin{array}{cc}
 A & B \\ C & D
\end{array}
\right)$.
Then $A,B,C,D$ are holomorphic functions that are
locally defined on a coordinate neighborhood of $M^2$.
We have only to show that
$\varphi\circ \nu$ is real analytic 
at the point $p$ where $|h(p)|=1$,
since $S_f=\{p\in M^2\,;\, |h(p)|=1\}$.  

By a straightforward calculation
using \eqref{eq:herm-mink}, 
\eqref{eq:fuji}, \eqref{eq:fuji1} and \eqref{eq:phi}, 
one can check that
each component of
$\varphi\circ \nu$ is a rational expression 
in $A,B,C,D,h$ and their conjugates,   
whose denominator $r$ is the following: 
\begin{eqnarray*}
r=
2(1-|h|^2)+|A+B\bar h|^2+|C+D\bar h|^2 
+|Ah+B|^2+|Ch+D|^2. 
\end{eqnarray*}
Hence, it is sufficient to show that
$r(p) \ne 0$ when $|h(p)|=1$.
However, the equality
\[
 Ah+B=A+B\bar h=C+D\bar h=Ch+D=0
\]
holds at $p$  only when $F(p)$ is a singular
matrix. 
Therefore $r(p)$ cannot be zero.  
This proves the 
assertion.
\end{proof}

\begin{corollary}\label{cor:smooth}
Let $f : M^2\to S^3_1$ be a CMC-1 face.
Then there exists a real analytic 
normal vector field defined on $M^2$.
In particular, $f$ is a 
co-orientable frontal.
\end{corollary}
\begin{proof}
Define a map $\psi: H^3_+ \cup H^3_- \to \R^4_1$ by 
\[
 \psi(u) = \begin{cases}
	   u /|u|_E & \text{if $u \in H^3_+$}, \\ 
	   -u /|u|_E & \text{if $u \in H^3_-$}, 
	  \end{cases}
\]
where $| u |_E$ denotes the Euclidean norm of  $u$ as 
a vector in $\R^4(=\R^4_1)$. 
For the unit normal field 
$\nu : M^2 \setminus S_f \to H^3_+ \cup H^3_- $, 
the equality $\psi \circ \nu = \psi \circ \varphi^{-1} \circ N$ 
holds on $M^2 \setminus S_f$. 
Here, $N : M^2 \to \bar{\mathcal H}^3$ denotes 
the real analytic extension of $\varphi \circ \nu : 
M^2 \setminus S_f \to \bar{\mathcal H}^3$ explained in Theorem 
\ref{thm:smooth}. %
On the other hand, we can calculate 
$\psi\circ \varphi^{-1} : \bar{\mathcal H}^3\setminus \{|x|=1\} 
\to \R^4_1$ 
as follows: 
\[
 \psi\circ \varphi^{-1}: (x,y,z)\mapsto 
\frac{1}{\sqrt{\Delta}}(1+x^2+y^2+z^2,-2x,-2y,-2z)
\]
where
\[
\Delta:=
\left(x^2+y^2+z^2+1\right)^2+4 x^2+4 y^2+4 z^2.
\]
It follows that $\psi\circ \varphi^{-1}$ can be smoothly extended to
the whole of $\bar{\mathcal H}^3$. Therefore 
$\Psi := \psi \circ \varphi^{-1} \circ N$ 
can be considered as a smooth map on $M^2$ such that 
$\Psi(p) \in \R^4_1$ 
is perpendicular to the position vector $f(p)\in \R^4_1$
and the subspace $df(T_pM)(\subset \R^4_1)$
at every regular point 
$p \in M^2 \setminus S_f$. 
Thus $\Psi$ gives a globally defined normal vector field 
of $f$,  which 
proves the co-orientability of $f$.  
\end{proof}

\section{Co-orientability of maximal surfaces in 
$\R^3_1$}\label{section:maxface}

In this section, we shall discuss 
the co-orientability
of maximal surfaces in Lorentz-Minkowski 3-space
given as a joint work with Shoichi Fujimori and Kotaro Yamada
during the authors' stay at RIMS in Kyoto 
and a meeting of the second author with Kotaro Yamada
at Osaka University on June 2009:
CMC-1 faces in $S^3_1$ have very similar properties
to maximal surfaces in Lorentz-Minkowski 3-space
$\R^3_1$ of signature $(-,+,+)$.
In \cite{UY2}, maxfaces (spacelike maximal surfaces
with admissible singular points) 
are defined as $C^\infty$ maps of 
orientable 2-manifolds into  $\R^3_1$.
If $M^2$ is non-orientable, we can take a double
covering $\hat \pi:\hat M^2\to M^2$ such that $\hat M^2$
is orientable.
For a non-orientable 2-manifold $M^2$,
a $C^\infty$ map $f:M^2\to \R^3_1$ is 
called a {\it maxface} if $f\circ \pi$ is a maxface
in the sense of \cite{UY2}.

\begin{figure}[ht]
\begin{center}
        \includegraphics[height=3.3cm]{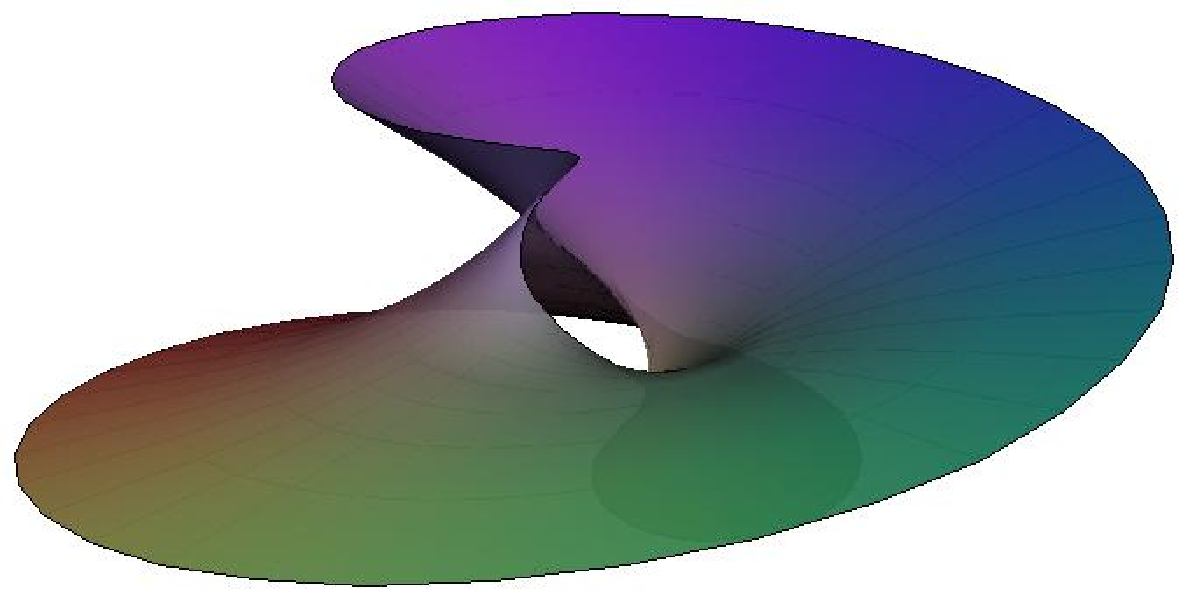}\hspace{1cm}
        \includegraphics[height=4.0cm]{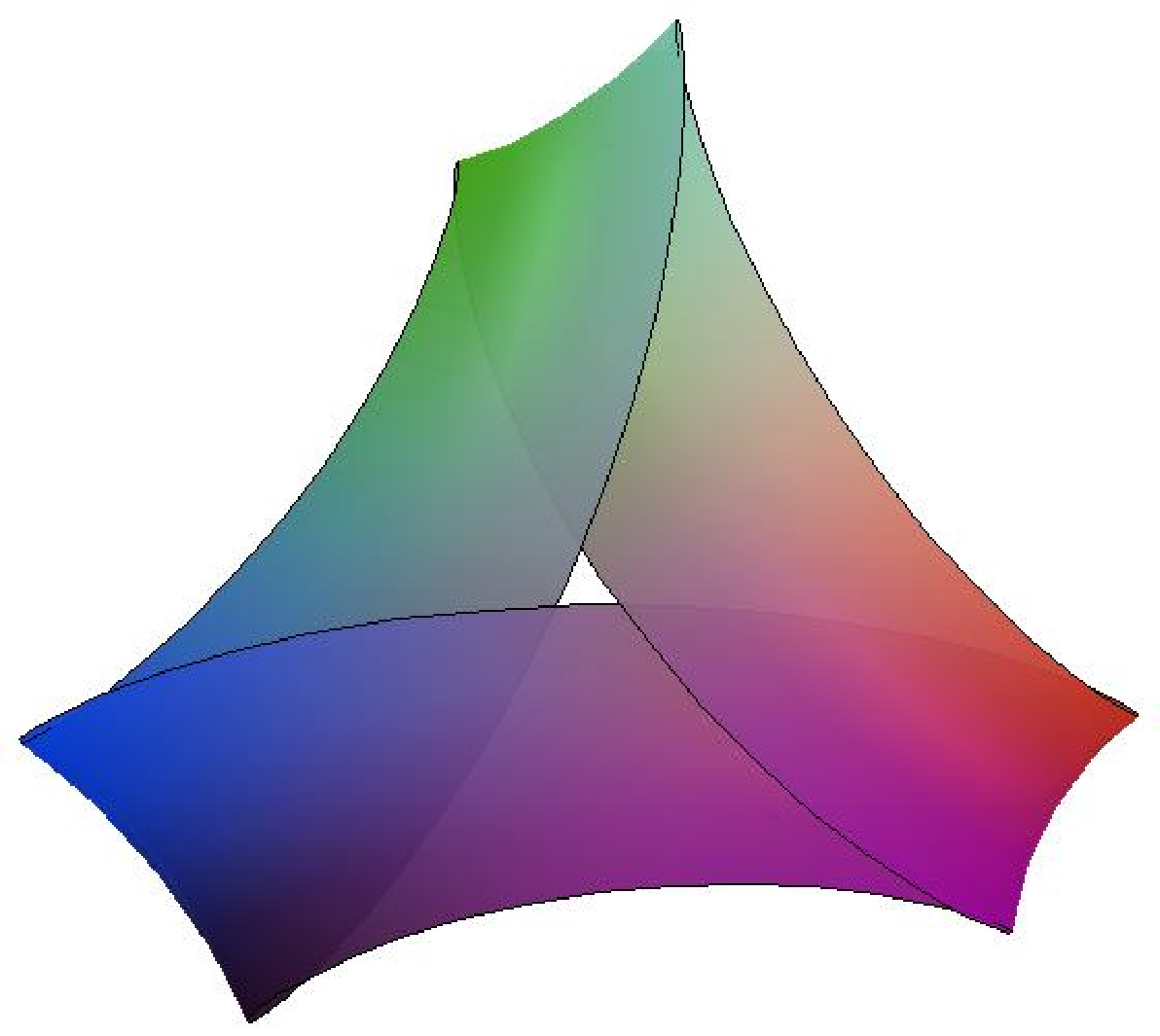}
\caption{Non-orientable maxfaces given in \cite{FL}
({the picture is produced by Fujimori})
}\label{non-ori}
\end{center}
\end{figure}

In contrast to CMC-1 faces in $S^3_1$,
there are non-orientable complete maxfaces
with or 
without handles (see 
Fujimori and L\'opez \cite{FL}).
Figure \ref{non-ori} shows
two different non-orientable maxfaces
given in \cite{FL}.
However, orientability and 
co-orientability of maxfaces are not discussed explicitly 
in \cite{UY2}.

The {\it hyperbolic $2$-sphere} $\bar {\mathcal H}^2$ is 
a Riemann sphere 
\[
 \bar {\mathcal H}^2:=\C\cup\{\infty\} \cong S^2
\]
endowed with the metric
$4 |dw|^2/(1-|w|^2)^2$ 
on $S^2 \setminus \{ {\rm equator} \}$, 
where $w \in \C \cup \{ \infty \}$. 
The hyperboloid of two sheets 
in $\R^3_1$, denoted by $H^2_+ \cup H^2_{-}$, is  
\[
H^2_+\cup H^2_- := \left\{ (x_0,x_1,x_2) \in \R^3_1 \, ; \, 
-x_0^2+x_1^2+x_2^2=-1 \right\}.
\] 
We consider the stereographic projection
\begin{equation*}
 \varphi:H^2_+\cup H^2_- \ni (x_0,x_1,x_2)\mapsto 
\frac{x_1+ix_2}{1-x_0}\in 
\bar {\mathcal H}^2,
\end{equation*}
which is an isometric embedding, 
and $\bar {\mathcal H}^2$ can be considered
as a compactification of $H^2_+\cup H^2_-$.
The following fact is well-known (cf. \cite{UY2}):
\begin{fact}\label{thm:smooth2}
Let $M^2$ be a Riemann surface,  
$f:M^2\to \R^3_1$  a maxface,
and $G:M^2\to \bar {\mathcal H}^2 =\C\cup\{\infty\}$
its Lorentzian Gauss map as a meromorphic function.
Then $\nu=\varphi^{-1}\circ G$ gives a unit
normal vector field defined on the regular
set of $f$.
\end{fact}
Using this, we shall prove the following assertion 
regardless of whether $M^2$ is orientable or not.

\begin{proposition}\label{thm:FKUY} {\rm (\cite{FKUY})}
Let $f:M^2\to \R^3_1$ be a maxface. 
Then there exists a real analytic 
normal vector field defined on $M^2$.
In particular, $f$ is a 
co-orientable frontal.
\end{proposition}

\begin{proof}
A smooth map $F:M^2\to \R^3_1$ is a co-orientable
frontal if and only if there exists a 
$C^\infty$ vector field $\eta$ defined on $M^2$
such that $\eta_p\in T_{F(p)}\R^3_1$ 
is Lorentzian orthogonal to the vector space $dF(T_pM^2)$.

We define a map $\psi: H^2_+ \cup H^2_- \to \R^3_1$ by 
\[
 \psi(u) = \begin{cases}
	   u /|u|_E & \text{if $u \in H^2_+$}, \\ 
	   -u /|u|_E & \text{if $u \in H^2_-$}, 
	  \end{cases}
\]
where $| u |_E$ denotes the Euclidean norm of  $u$ as 
a vector in $\R^3(=\R^3_1)$. 

Suppose that $M^2$ is orientable.
By Fact \ref{thm:smooth2},
\begin{equation}\label{eq:hat_nu}
\nu:=\psi\circ \varphi^{-1}\circ G
\end{equation}
gives a Lorentzian normal vector field of $f$
on $M^2 \setminus S_{f}$. 
On the other hand, we can calculate 
$\psi\circ \varphi^{-1} : \bar{\mathcal H}^2\setminus \{ |w|=1\} 
\to \R^3_1$ as follows: 
\[
 \psi\circ \varphi^{-1}: w=x+iy \mapsto 
\frac{1}{\sqrt{\Delta}}(1+x^2+y^2,-2x,-2y),
\]
where
\[
\Delta:=
\left(x^2+y^2+1\right)^2+4 x^2+4 y^2.
\]
It follows that $\psi\circ \varphi^{-1}$ can be smoothly extended 
to the whole of $\bar{\mathcal H}^2$. Therefore 
$\nu$ given in \eqref{eq:hat_nu}
is a globally defined normal vector field 
of $f$. 
In fact, by a straightforward calculation, we have
that
\begin{equation}\label{eq:N}
\nu=\frac{1}{\sqrt{(1+|G|^2)^2+4|G|^2}}
\left(1+|G|^2,-2\operatorname{Re}(G),-2\operatorname{Im}(G)
\right).
\end{equation}
So the orientability implies 
the co-orientability. 
Next, suppose that  $f:M^2\to \R^3_1$
is a maxface and $M^2$ is non-orientable.
Then there exists a double covering 
$\hat\pi : \hat M^2 \to M^2$ such that
$\hat M^2$ is orientable.
Set $\hat f:=f\circ \hat \pi$
and denote by $G:M^2\to \C\cup\{\infty\}$
the Gauss map of the maxface $\hat f$.
Let $T:\hat M^2\to \hat M^2$ be the covering involution
of the double covering $\hat \pi:\hat M^2\to M^2$.
Then $(\varphi^{-1} \circ G) \circ T(p)$ is antipodal to 
$(\varphi^{-1} \circ G)(p)$ in 
$H^2_+ \cup H^2_-$ for every $p \in \hat M^2$ (cf. \cite{FL}). 
It follows that 
\begin{equation}\label{eq:GT}
G\circ T(p)={1} / \, {\overline{G(p)}}\qquad (p\in \hat M^2).
\end{equation}
By \eqref{eq:GT}, 
the vector field $\nu$ given by \eqref{eq:N}
is a Lorentzian normal vector field
on $\hat M^2$
satisfying $\nu\circ T=\nu$.
In other words, $\nu$ can be considered as a smooth 
normal vector field on $M^2$, which
proves the assertion.
\end{proof}

\begin{remark}
Let $f:M^2\to \R^3_1$ be a non-orientable maxface.
Then a closed regular curve
$
\gamma:[0,1]\to M^2
$
is called a {\it non-orientable loop} if
the orientation of $M^2$ reverses along the curve.
Then as pointed out in Fujimori and L\'opez \cite{FL},
there is at least one singular point on $\gamma$.
Take a double covering 
$\hat\pi :\hat M^2 \to M^2$ such that
$\hat M^2$ is orientable.
Set $\hat f:=f\circ \hat \pi$
and denote by $G:\hat M^2\to \bar {\mathcal H}^2=\C\cup\{\infty\}$
the Lorentzian Gauss map of the maxface $\hat f$.
Since the singular set $S_{\hat f}$ of $\hat f$ is the $1$-level set
of $|G|$,
it consists of a union of regular curves outside of the
discrete set
\[
\Gamma:=\left\{p\in S_{\hat f}\,;\, dG(p)=0\right\}.
\]
Here $G$ is considered as a smooth map into $S^2$
and $dG$ is the derivative of $G$.
A non-orientable loop $\gamma$ 
is called {\it generic} if its lift $\hat \gamma$ on $\hat M^2$
meets the singular set only at
finitely many points on $S_{\hat f}\setminus \Gamma$,
where $\hat \gamma$ passes through these points 
transversally to the set $S_{\hat f}\setminus \Gamma$.
We fix a generic non-orientable loop $\gamma$ arbitrarily.
Since $T\circ \hat \gamma (0)=\hat \gamma(1)$,
by \eqref{eq:GT}, we have 
$|G\circ \hat \gamma (0)| = 
1/ |G \circ \hat \gamma (1)|$.
Since the curve $G\circ \hat \gamma(t)$ 
$(t\in [0,1])$ is continuous,
it meets the 
circle $S^1 = \{ w \in \C \cup \{ \infty \} \, ; |w|=1 \}$ 
odd times.
Thus {\it the number of singular points
on a generic non-orientable loop $\gamma$ is
odd.}
(In fact, the number of connected components of the
singular set on the non-orientable maxfaces
in Figure \ref{non-ori} is one and three, 
respectively.)
\end{remark}

\appendix

\section{
Degenerate Weingarten fronts}
In hyperbolic $3$-space $H^3$, the following
facts are well-known;
\begin{enumerate}
 \item[(1)] one of the parallel surfaces of an open portion of 
a hyperbolic cylinder 
degenerates to a hyperbolic line, and  
 \item[(2)] one of the parallel surfaces of an open portion of 
a geodesic sphere 
(a compact totally umbilical surface) degenerates to a 
single point. 
\end{enumerate}
Note that these are Weingarten fronts. 
(A {\it Weingarten front} is defined by Definition \ref{def:2}  
removing the adjective `horospherical linear'. 
This definition makes sense based on the fact that any 
parallel surface of a Weingarten surface is again Weingarten.)  
The above fronts (1) and (2) are called 
a {\it degenerate cylinder} and 
a {\it degenerate sphere\/}, respectively. 
We wish to prove a converse statement. Before this, 
let us recall the following well known lemma: 
\begin{lemma}\label{lem:parallel-sing}
 Let $f : M^2 \to H^3$ be an immersion, and 
$\kappa_i$ $(i=1,2)$ the principal curvatures. 
Then the set of singular points  
of a parallel surface $f_\delta$ coincides with   
$\{ p \in M^2 \, ; \, \coth^{-1} \kappa_i(p) = \delta \}$.  
\end{lemma} 
Lemma \ref{lem:parallel-sing} asserts that $p \in M^2$ 
is a singular point of $f_{\delta}$ 
if and only if $\delta$ is the radius $\coth^{-1} \kappa_i(p)$ 
of the principal curvature $\kappa_i(p)$.  
It is easily observed that
 any parallel surface is free from singular points if 
and only if $|\kappa_i| \le 1$.  

\begin{proposition}
Let $f:M^2 \to H^3$ be a Weingarten front. 
Then all points in $M^2$ are singular points if and only if 
$f$ is locally a degenerate cylinder or a degenerate sphere. 
\end{proposition}
\begin{proof}
Suppose that all points in $M^2$ are singular points of $f$. 
Let $p \in M^2$ be an arbitrary point.  
There exist a neighborhood $U$ of $p$ and a real number 
$\delta_0$ such that $f_{\delta_0} : U \to H^3$ is an 
immersion. In other words, $f: U \to H^3$ is a 
parallel surface of an immersion $g:=f_{\delta_0}$, i.e., 
$f=g_{-\delta_0}$.  
Since all points are singular, 
it follows from Lemma \ref{lem:parallel-sing}
that  
\begin{equation*}
  -\delta_0 = \coth^{-1} \kappa(p) \ \text{ for all }p \in U, 
\end{equation*}  
where $\kappa$ denotes one of the principal curvatures of $g$. 
It implies that $\kappa$ is constant on $U$ and satisfies 
$|\kappa|>1$. Moreover, another principal 
curvature is also constant because $g$ is Weingarten. 
Thus $g$ is an isoparametric surface. By the classification 
due to Cartan \cite{cartan}, the immersion $g$ is totally umbilical or 
locally a standard embedding of $S^1 \times H^1$ into $H^3$. 
Under our assumption, the former is a part of geodesic sphere because 
$|\kappa|>1$, the latter is a part of a hyperbolic cylinder. 

Finally, by the connectivity of $M^2$ and by the continuity of $f$, 
we can conclude that $f:M^2 \to H^3$ is a degenerate cylinder or 
a degenerate sphere.    
\end{proof}

\medskip
\noindent
{\bf Acknowledgement}
The authors thank Antonio Mart\'\i{}nez
for a fruitful discussion on 
the geometry of linear Weingarten surfaces, when they visited
the Fukuoka conference on
November, 2004.
They also thank Shyuichi Izumiya
for his informal lecture
on  the duality
 between fronts in $H^3$ and $S^3_1$
at the Karatsu workshop on
October, 2005.
The first author thanks Shoichi Fujimori for 
fruitful discussions about the proof of 
Proposition \ref{prop:Gstar},  
after the Karatsu workshop.
The second author thanks Masahiko Kanai, who
pointed out the importance of $L(H^3)$ 
as the range of hyperbolic Gauss maps of CMC-1 surfaces
in $H^3$ when the second author stayed  
at Nagoya University in 2001. 
Finally, the authors thank Kotaro Yamada, Wayne Rossman,
 Francisco L\'opez and the referee for valuable comments.


\begin{thebibliography}{99}
\bibitem{AE}
 J.~A.~Aledo and J.~M.~Espinar,
  \textrm{A conformal representation for linear Weingarten surfaces
in the de Sitter space},
J. Geom. Phys. 
\textbf{57}, 1669--1677 (2007).
%
\bibitem{AG}
 J.~A.~Aledo and J.~A.~G\'alvez,
  \textrm{Complete surfaces in hyperbolic space with a 
constant principal curvature},
Math. Nachr.   
\textbf{278}, 1111--1116 (2005).
%
\bibitem{Br}
 R.~Bryant,
  \textrm{Surfaces of mean curvature one in hyperbolic space},
  Ast\'erisque \textbf{154--155}, 321--347 (1987).
%
\bibitem{cartan}
\'E. Cartan, 
 \textrm{Familles de surfaces isoparam\'etriques dans les espaces 
\`a courbure constante}, 
Ann. Mat. Pura Appl. \textbf{17}, 177--191 (1938). 
%
\bibitem{chen}
  B. Y. Chen,
  \textrm{Minimal surfaces with
         constant Gauss curvature}, 
  Proc. Amer. Math. Soc. \textbf{34}, 504--508 (1972).
%
\bibitem{CFG}
V.~Cruceanu, P.~Fortuney and P.~M.~Gadea, 
\textrm{A survey on paracomplex geometry},  
Rocky Mountain J. Math. \textbf{26}, 83--115 (1996).
%
\bibitem{F}
  S. Fujimori,
  \textrm{Spacelike CMC $1$ surfaces with
         elliptic ends in de Sitter $3$-Space},  
  Hokkaido Math. J. \textbf{35}, 289--320 (2006).
%
\bibitem{FL}
S.~Fujimori \and F.~J.~L\'opez, 
Nonorientable maximal surfaces in the 
Lorentz-Minkowski $3$-space, 
Tohoku Math. J. \textbf{62}, 311--328 (2010). 
%
\bibitem{FRUYY}
S.~Fujimori, W.~Rossman, M.~Umehara, K.~Yamada \and S.-D.~Yang, 
        \textrm{Spacelike mean curvature one surfaces in 
de Sitter $3$-space}, Comm. Anal. Geom. \textbf{17}, 
383--427 (2009).
%
\bibitem{FSUY}
  S. Fujimori, K. Saji, M. Umehara \and K. Yamada,
\textrm{Singularities of maximal surfaces}, 
Math. Z. \textbf{259}, 827--848 (2008).
%
\bibitem{FKUY}
S. Fujimori, M. Kokubu,  M. Umehara and K. Yamada,
\textrm{Personal meetings at RIMS Kyoto
and at Osaka University on June, 2009}. 
%
\bibitem{GG}
N.~Georgiou \and B.~Guilfoyle,
\textrm{A characterization of Weingarten surfaces in
hyperbolic $3$-space}, 
Abh. Math. Semin. Univ. Hambg. \textbf{80}, 233-253 (2010).
%
\bibitem{GMM1}
 J. A. G\'alvez, A. Mart\'\i{}nez \and F. Mil\'an,
  \textrm{Flat surfaces in hyperbolic $3$-space},  
  Math. Ann. \textbf{316}, 419--435 (2000).
%
\bibitem{GMM2}
J. A. G\'alvez, A. Mart\'\i{}nez \and F. Mil\'an, 
\textrm{Complete linear Weingarten surfaces of Bryant type. 
A Plateau problem at infinity},  
Trans. Amer. Math. Soc. \textbf{356}, 3405--3428 (2004).
%
\bibitem{IPRT} 
S. Izumiya, D. Pei, M. C. Romero-Fuster 
\and M. Takahashi,
\textrm{The horospherical geometry of submanifolds in
hyperbolic space}, 
J. London Math. Soc. (2)  \textbf{71}, 779--800  (2005).
%
\bibitem{IS} 
S. Izumiya and K. Saji, 
\textrm{The mandala of Legendrian dualities for pseudo-
spheres in Lorentz-Minkowski space and ``flat'' spacelike surfaces}, 
J. Singul. \textbf{2}, 92--127  (2010).

\bibitem{IST} 
S. Izumiya, K. Saji \and M. Takahashi, 
\textrm{Horospherical flat surfaces in Hyperbolic 
$3$-space}, 
J. Math. Soc. Japan \textbf{62}, 789--849 (2010).
%
\bibitem{Kanai} 
M. Kanai,
\textrm{Geodesic flows of negatively curved manifolds with 
smooth stable and unstable foliations},  
Ergodic Theory Dynam. Systems \textbf{8}, 215--239 (1988). 
%
\bibitem{Kane} 
S. Kaneyuki and M. Kozai,
\textrm{Paracomplex structures and affine symmetric spaces},  
Tokyo J. Math. \textbf{8}, 81--98 (1985). 
%
\bibitem{Kok} 
M.~Kokubu, 
\textrm{Surfaces and fronts 
with harmonic-mean curvature one in hyperbolic three-space},
Tokyo J. Math. \textbf{32}, 177--200 (2009).
%
\bibitem{KRSUY}
M.~Kokubu, W.~Rossman, K.~Saji, M.~Umehara \and K.~Yamada, 
\textrm{Singularities of flat fronts in hyperbolic space},  
Pacific J. Math. \textbf{221}, 303--351 (2005).
%
\bibitem{KRUY}
 M. Kokubu, W. Rossman, M. Umehara \and K. Yamada,  
\textrm{Flat fronts in hyperbolic $3$-space and their caustics}, 
 J. Math. Soc. Japan \textbf{59}, 265--299 (2007). 
%
\bibitem{KUY}
        M.\ Kokubu, M.\ Umehara and K.\ Yamada,
        \textrm{Flat fronts  in hyperbolic $3$-space},
        Pacific J.\ Math. {\bf 216}, 149--175 (2004).
%
\bibitem{LLR}
R.\ Langevin, G.\ Levitt \and H.\ Rosenberg,
 \textrm{Classes d'homotopie de surfaces avec rebroussements
    et queues d'aronde dans $\boldsymbol R^3$}, 
    Canad.\ J. Math. \textbf{47}, 544--572 (1995).

\bibitem{SUY1}
K. Saji, M. Umehara \and K. Yamada,
 \textrm{The geometry of fronts},
Ann.\ of Math. {\bfseries 169}, 491--529 (2009). 
%
\bibitem{SUY2}
K. Saji, M. Umehara \and K. Yamada,
\textrm{$A_k$ singularities of wave fronts}, 
Math. Proc. Cambridge Philos. Soc. \textbf{146}, 731--746 (2009).
\bibitem{TT} 
C. Takizawa \and K. Tsukada, 
\textrm{Horocyclic surfaces in Hyperbolic $3$-space}, 
Kyushu J. Math. \textbf{63}, 269--284 (2009).
%
\bibitem{UY1}
M.~Umehara \and K.~Yamada,
  \textrm{Complete surfaces of constant mean curvature-$1$
         in the hyperbolic $3$-space},  
  Ann.\ of Math. \textbf{137}, 611--638 (1993).


\bibitem{UY2}
M.~Umehara \and K.~Yamada,
 \textrm{Maximal surfaces with singularities in
     Minkowski space}, Hokkaido Math. J. 
\textbf{35}, 13--40 (2006).
\end{thebibliography}
\end{document}